\documentclass[12pt]{interact}

\usepackage{epstopdf}% To incorporate .eps illustrations using PDFLaTeX, etc.
\usepackage[caption=false]{subfig}% Support for small, `sub' figures and tables
\usepackage[doublespacing]{setspace}% To produce a `double spaced' document if required
\usepackage{amsmath}
\usepackage{color}

\usepackage[numbers,sort&compress]{natbib}% Citation support using natbib.sty
\bibpunct[, ]{[}{]}{,}{n}{,}{,}% Citation support using natbib.sty
\makeatletter% @ becomes a letter
\def\NAT@def@citea{\def\@citea{\NAT@separator}}% Suppress spaces between citations using natbib.sty
\makeatother% @ becomes a symbol again

\theoremstyle{plain}% Theorem-like structures provided by amsthm.sty

\theoremstyle{definition}

\theoremstyle{remark}

\begin{document}

%\articletype{ARTICLE TEMPLATE}% Specify the article type or omit as appropriate

%%%%%%%%%%%%%%%%%%%%%%%%%%%%%%%%%%%%%%%%%%%%%%%%%%%%%%%%%%%%%%%%%%%%%%%%%%%%%%%%
\title{A Bivariate Beta Distribution with Arbitrary Beta Marginals and its Generalization to a Correlated Dirichlet Distribution}
%%%%%%%%%%%%%%%%%%%%%%%%%%%%%%%%%%%%%%%%%%%%%%%%%%%%%%%%%%%%%%%%%%%%%%%%%%%%%%%%

\author{
\name{Susanne Trick\textsuperscript{a,b}\thanks{CONTACT Susanne Trick. Email: susanne.trick@cogsci.tu-darmstadt.de}, Frank J\"akel\textsuperscript{a,b,*}, and Constantin A. Rothkopf\textsuperscript{a,b,c,*}}
\affil{\textsuperscript{a}Centre for Cognitive Science, TU Darmstadt, Darmstadt, Germany; \textsuperscript{b}Institute of Psychology, TU Darmstadt, Darmstadt, Germany; \textsuperscript{c}Frankfurt Institute for Advanced Studies, Goethe University Frankfurt, Frankfurt, Germany}
\affil{\textsuperscript{*}Frank Jäkel and Constantin A. Rothkopf contributed equally.}
}

\maketitle

\begin{abstract}
We discuss a bivariate beta distribution that can model arbitrary beta-distributed marginals with a positive correlation. The distribution is constructed from six independent gamma-distributed random variates. We show how the parameters of the distribution can be fit to data using moment matching. Previous work used an approximate and sometimes inaccurate method to compute the covariance. Here, we derive all product moments and the exact covariance, which can easily be computed numerically. The bivariate case can be generalized to a multivariate distribution with arbitrary beta-distributed marginals. Furthermore, we generalize the distribution from two marginal beta to two marginal Dirichlet distributions. The resulting correlated Dirichlet distribution makes it possible to model two correlated Dirichlet-distributed random vectors.
\end{abstract}


%%%%%%%%%%%%%%%%%%%%%%%%%%%%%%%%%%%%%%%%%%%%%%%%%%%%%%%%%%%%%%%%%%%%%%%%%%%%%%%%
\section{Introduction}
%%%%%%%%%%%%%%%%%%%%%%%%%%%%%%%%%%%%%%%%%%%%%%%%%%%%%%%%%%%%%%%%%%%%%%%%%%%%%%%%

Many bivariate beta distributions have been proposed in the literature \cite{libby1982, nadarajah2005, jayakumar2019, sarabia2006, orozco2012, brancardona2011, gupta2011, gupta1985, tinglee1996, nadarajah2017, jones2002, elbassiouny2009, olkin2003, arnold2011, olkin2015, arnold2017, magnussen2004}. Among them, some approaches have been derived from general families of bivariate distributions, such as the Farlie-Gumbel-Morgenstern family of distributions (e.g. \cite{gupta1985}) and the Sarmanov family of distributions (e.g. \cite{tinglee1996}). Others derived a bivariate beta distribution from a bivariate extension of the F distribution \cite{jones2002, elbassiouny2009}, whereas Nadarajah and Kotz proposed different bivariate beta distributions constructed from products of univariate beta-distributed random variables \cite{nadarajah2005}. In general, using different copulas one can construct different bivariate distributions with the same beta marginals. However, as beta-distributed random variables can easily be constructed from normalized gamma-distributed random variables, it is natural to try and generalize this construction to the bivariate case. In this vein, several authors have introduced correlations through shared gamma-distributed random variables \cite{olkin2003, arnold2011, olkin2015, arnold2017, magnussen2004}. This is also the approach we will follow in this paper where we are interested in bivariate distributions with arbitrary beta-distributed marginals and a positive correlation.

The most straightforward case of this construction has been studied by Olkin and Liu \cite{olkin2003}, building on the work of Libby and Novick \cite{libby1982}. They use three independent gamma-distributed random variables $U_1, U_2, U_3$ with respective shape parameters $\upsilon_1, \upsilon_2, \upsilon_3$ and same scale parameter to construct
\begin{align}\label{eq:olkin_liu}
X' = \frac{U_1}{U_1+U_3} \quad\mathrm{and}\quad Y' = \frac{U_2}{U_2+U_3}.
\end{align}
Using the standard construction of beta variates from gamma variates, the joint distribution of the random variables $X'$ and $Y'$ is a bivariate beta distribution with marginal distributions Beta($X'$; $\upsilon_1$, $\upsilon_3$) and Beta($Y'$; $\upsilon_2$, $\upsilon_3$). The correlation between $X'$ and $Y'$, which is obtained through the shared latent variable $U_3$ and its parameter $\upsilon_3$, is in the range [0,1]. For high values of $\upsilon_3$, the correlation tends to 0 whereas for low values of $\upsilon_3$, it tends to 1. However, if $\upsilon_3$ is high, the values of $X'$ and $Y'$ also tend to 0 and if $\upsilon_3$ is low they tend to 1 accordingly. This behavior severely limits the usefulness of the distribution for most applications. A further limitation is the constraint that the marginal distributions share the same second parameter $\upsilon_3$. Thus, the bivariate beta distribution proposed by Olkin and Liu does not allow for arbitrary beta marginals.

Arnold et al. \cite{arnold2011} proposed a more flexible construction for a bivariate beta distribution. They use five independent gamma-distributed random variables $U_1, \dots, U_5$ with shape parameters $\upsilon_1, \dots, \upsilon_5$ and scale parameter 1 to define two correlated random variables
\begin{align}
X = \frac{U_1 + U_3}{U_1 + U_3 + U_4 + U_5} \quad\mathrm{and}\quad Y = \frac{U_2 + U_4}{U_2 + U_3 + U_4 + U_5}
\end{align}
with marginal distributions Beta($X$; $\upsilon_1 + \upsilon_3, \upsilon_4 + \upsilon_5$) and Beta($Y$; $\upsilon_2 + \upsilon_4, \upsilon_3 + \upsilon_5$). Compared to Olkin et al. \cite{olkin2003}, this construction of a bivariate beta distribution can generate all correlations in the range [-1,1] and marginal distributions with differing second parameters. Nevertheless, because of how the two marginals share parameters, not all combinations of parameters of the marginal beta distributions are possible. For example, the marginals Beta($X$; 10, 4) and Beta($Y$; 1, 1) cannot be obtained.

Olkin and Trikalinos \cite{olkin2015} base their construction of a bivariate beta distribution on the Dirichlet distribution. $U = (U_{00}, U_{01}, U_{10}, U_{11})$ is drawn from a 4-dimensional Dirichlet distribution with parameters $\upsilon_1, \dots, \upsilon_4$. By just using three of its components, $U_{00}, U_{01}, U_{10}$, new random variables
\begin{align}\label{eq:olkin_dirichlet}
X = U_{00} + U_{10} \quad\mathrm{and}\quad Y = U_{01} + U_{10}
\end{align}
are constructed, with marginal distributions Beta($X$; $\upsilon_1 + \upsilon_3, \upsilon_2 + \upsilon_4$) and Beta($Y$; $\upsilon_2 + \upsilon_3, \upsilon_1 + \upsilon_4$). As Dirichlet-distributed random variables can also be constructed from gamma random variables, we can equivalently construct $X$ and $Y$ in equation (\ref{eq:olkin_dirichlet}) from four independent gamma-distributed random variables $U_1, \dots, U_4$ with shape parameters $\upsilon_1, \dots, \upsilon_4$ and equal scale parameter 1, with
\begin{align} \label{eq:olkin_gamma}
X = \frac{U_1 + U_3}{U_1 + U_2 + U_3 + U_4} \quad\mathrm{and}\quad Y = \frac{U_2 + U_3}{U_1 + U_2 + U_3 + U_4}.
\end{align}
As can easily be seen from this construction, all correlations in the range [-1,1] can be generated. In particular, the correlation tends to -1 if $\upsilon_3$ and $\upsilon_4$ tend to 0 and $U_3$ and $U_4$ will be negligible compared to $U_1$ and $U_2$. In this case $X \approx \frac{U_1}{U_1 + U_2} \approx 1 - Y$. Similarly, the higher the values of $\upsilon_3$ and $\upsilon_4$ relative to $\upsilon_1$ and $\upsilon_2$, the more negligible $U_1$ and $U_2$ will be and the correlation increases to 1 until $X \approx \frac{U_3}{U_3 + U_4} \approx Y$. Less obviously, a correlation of 0 is obtained in case $\upsilon_1 \cdot \upsilon_2 = \upsilon_3 \cdot \upsilon_4$ \cite{olkin2015}. Still, this construction of a bivariate beta distribution does not allow arbitrary beta marginal distributions. Since all $\upsilon_i$ are constrained to be positive, for some combinations of marginal distributions the resulting system of linear equations for the parameters $\upsilon_i$ has no solution. For example, the two marginals Beta($X$; 2, 2) and Beta($Y$; 1, 1) cannot be generated, regardless of their correlation.

Magnussen \cite{magnussen2004} introduced yet another construction based on six gamma variates. While all the constructions thus far constrain the parameters of the marginal beta distributions, this construction does allow for arbitrary beta marginals. But it only allows modeling positive correlations. Magnussen's distribution \cite{magnussen2004} is a special case of a more general 8-parameter bivariate beta distribution introduced by Arnold and Ng \cite{arnold2011} and reviewed in Arnold and Gosh \cite{arnold2017}. As shown in this review, the 8-parameter distribution subsumes all the other constructions we have discussed above \cite{olkin2003, arnold2011, olkin2015} as special cases and allows for positive and negative correlations. However, in many applications it is enough to model positive correlations, for which the less complex 6-parameter distribution is sufficient. For example, if $X$ and $Y$ are probability estimates elicited from two skilled forecasters, we do not expect negative correlations. But we do want to allow for the possibility that their marginal forecasts have different distributions that should not be tied together by parameter constraints on the marginals. Hence, in this paper we will focus on the bivariate beta distribution proposed by Magnussen \cite{magnussen2004} that can model arbitrary beta-distributed marginals with a positive correlation.

In the following, we will derive the full joint distribution, which has been missing in \cite{magnussen2004}, and we will clarify the relationship between Magnussen's distribution and the Olkin-Liu distribution \cite{olkin2003}. While previous work on the distribution only derived a rough and sometimes inaccurate approximation for the covariance \cite{magnussen2004}, using this link we derive all product moments and the exact covariance function (and we correct a small mistake in the product moments from \cite{olkin2003}). For parameter inference, we propose to match moments numerically using the exact covariance, because it is not available in closed-form. In order to do so we first show how to reasonably constrain the distribution’s parameters and thereby make parameter inference unambiguous.
We then extend the construction to the multivariate case, also with arbitrary marginal beta distributions. We conclude the paper by generalizing the distribution from two marginal beta to two marginal Dirichlet distributions. The resulting correlated Dirichlet distribution can model correlated Dirichlet-distributed random variables and is easy to interpret.

%%%%%%%%%%%%%%%%%%%%%%%%%%%%%%%%%%%%%%%%%%%%%%%%%%%%%%%%%%%%%%%%%%%%%%%%%%%%%%%%
\section{A Bivariate Beta Distribution with Arbitrary Beta Marginals} \label{sec:6_params}
%%%%%%%%%%%%%%%%%%%%%%%%%%%%%%%%%%%%%%%%%%%%%%%%%%%%%%%%%%%%%%%%%%%%%%%%%%%%%%%%

Magnussen \cite{magnussen2004} uses six independent gamma-distributed random variables $A_1, A_2, B_1, B_2, D_1, D_2$, which are distributed according to
\begin{align} \label{eq:biv_beta_6_variables}
\begin{split}
A_i &\sim \text{Gamma}(\alpha_i, 1)\quad i = 1,2 \\
B_i &\sim \text{Gamma}(\beta_i, 1)\hspace{0.02cm}\quad i = 1,2 \\
D_i &\sim \text{Gamma}(\delta_i, 1)\hspace{0.06cm}\quad i = 1,2, 
\end{split}
\end{align}
to construct two bivariate-beta-distributed random variables
\begin{align} \label{eq:bivariate_beta_6}
X &= \frac{A_1 + D_1}{A_1 + A_2 + D_1 + D_2} \quad\mathrm{and}\quad Y = \frac{B_1 + D_1}{B_1 + B_2 + D_1 + D_2}.
\end{align}
The resulting marginal distributions of $X$ and $Y$ are Beta($X; a_1, a_2$) and Beta($Y; b_1, b_2$) with
\begin{align}\label{eq:ab}
a_1 &= \alpha_1 + \delta_1, & a_2 &= \alpha_2 + \delta_2, & b_1 &= \beta_1 + \delta_1 & b_2 &= \beta_2 + \delta_2. 
\end{align}

The marginals follow immediately from the definition because the sum of gamma random variables of the same scale is gamma-distributed with the same scale but with the original shape parameters summed. In contrast to other constructions that were discussed above \cite{olkin2003, arnold2011,olkin2015}, this construction allows for arbitrary marginal distributions. In particular, when $\delta_1$ and $\delta_2$ tend to zero, we can model arbitrary independent marginal distributions Beta($X; \alpha_1, \alpha_2$) and Beta($Y; \beta_1, \beta_2$).

Since all parameters $\alpha_1, \alpha_2, \beta_1, \beta_2, \delta_1, \delta_2$ need to be positive by definition, for fixed marginal distributions Beta($X; a_1, a_2$) and Beta($Y; b_1, b_2$) it must hold that $\delta_1 < \text{min}(a_1, b_1)$ and $\delta_2 < \text{min}(a_2, b_2)$. Therefore, for most marginal distributions the maximum correlation that can be generated is below 1. The higher the difference between two marginal distributions, the lower the possible maximum correlation. A perfect correlation approaching 1 can, of course, only be generated for equal marginal distributions, i.e. if $a_1 = b_1$ and $a_2 = b_2$ and $\alpha_1$, $\alpha_2$, $\beta_1$, and $\beta_2$ tend to 0, as also noted by Magnussen \cite{magnussen2004}. However, note that this limitation applies to other bivariate distributions that do not allow for arbitrary marginal beta distributions as well (see e.g. \cite{olkin2015}).

The construction of this bivariate beta distribution can be seen as a pairwise combination of three beta distributions, Beta($\alpha_1, \alpha_2$), Beta($\beta_1, \beta_2$) and Beta($\delta_1, \delta_2$). First transform the six independent gamma-distributed random variables (\ref{eq:biv_beta_6_variables}) into three independent gamma- and three independent beta-distributed random variables,
\begin{equation}
\begin{aligned}\label{eq:uw}
U_{1} & =A_{1}+A_{2}, & U_{1} & \sim\text{Gamma}(\upsilon_{1},1)\\
U_{2} & =B_{1}+B_{2}, & U_{2} & \sim\text{Gamma}(\upsilon_{2},1)\\
U_{3} & =D_{1}+D_{2}, & U_{3} & \sim\text{Gamma}(\upsilon_{3},1)\\
W_{1} & =\frac{A_{1}}{A_{1}+A_{2}}, & W_{1} & \sim\text{Beta}(\alpha_{1},\alpha_{2})\\
W_{2} & =\frac{B_{1}}{B_{1}+B_{2}}, & W_{2} & \sim\text{Beta}(\beta_{1},\beta_{2})\\
W_{3} & =\frac{D_{1}}{D_{1}+D_{2}}, & W_{3} & \sim\text{Beta}(\delta_{1},\delta_{2})
\end{aligned}
\end{equation}
with
\begin{align}\label{eq:upsilon}
\upsilon_{1}&=\alpha_{1}+\alpha_{2}, & \upsilon_{2}&=\beta_{1}+\beta_{2} & \upsilon_{3}&=\delta_{1}+\delta_{2}.
\end{align}
With these definitions we can then rewrite construction (\ref{eq:bivariate_beta_6}) as
\begin{align} \label{eq:uvw}
\begin{split}
X & =\frac{U_{1}}{U_{1}+U_{3}}\cdot W_{1}+\frac{U_{3}}{U_{1}+U_{3}}\cdot W_{3}=X^{\prime}W_{1}+(1-X^{\prime})W_{3}\\
Y & =\frac{U_{2}}{U_{2}+U_{3}}\cdot W_{2}+\frac{U_{3}}{U_{2}+U_{3}}\cdot W_{3}=Y^{\prime}W_{2}+(1-Y^{\prime})W_{3},
\end{split}
\end{align}
where $X'$ and $Y'$ are defined as in (\ref{eq:olkin_liu}) but with $\upsilon_1$, $\upsilon_2$, and $\upsilon_3$ as in (\ref{eq:upsilon}). Furthermore, $X'$ and $Y'$ are independent of $W_1, W_2, W_3$. If parameters $\delta_1$ and $\delta_2$ and with them $U_3$ tend to 0, $X \approx W_1$ and $Y \approx W_2$ are independent with marginal distributions Beta($X; \alpha_1, \alpha_2$) and Beta($Y; \beta_1, \beta_2$). Mixing in the shared component $W_3$ by increasing the values of parameters $\delta_1$ and $\delta_2$ increases the correlation between $X$ and $Y$. If $U_1$ and $U_2$ are negligible compared to $U_3$ because $\delta_1$ and $\delta_2$ dominate the parameters, the correlation will be close to 1 with $X \approx W_3 \approx Y$ and hence $X$ and $Y$ have the same marginal distribution Beta($\delta_1, \delta_2$), as mentioned before.

%%%%%%%%%%%%%%%%%%%%%%%%%%%%%%%%%%%%%%%%
\subsection{Joint Distribution} \label{sec:6_params_joint}
%%%%%%%%%%%%%%%%%%%%%%%%%%%%%%%%%%%%%%%%

\begin{figure}
	\includegraphics[width=0.95\textwidth]{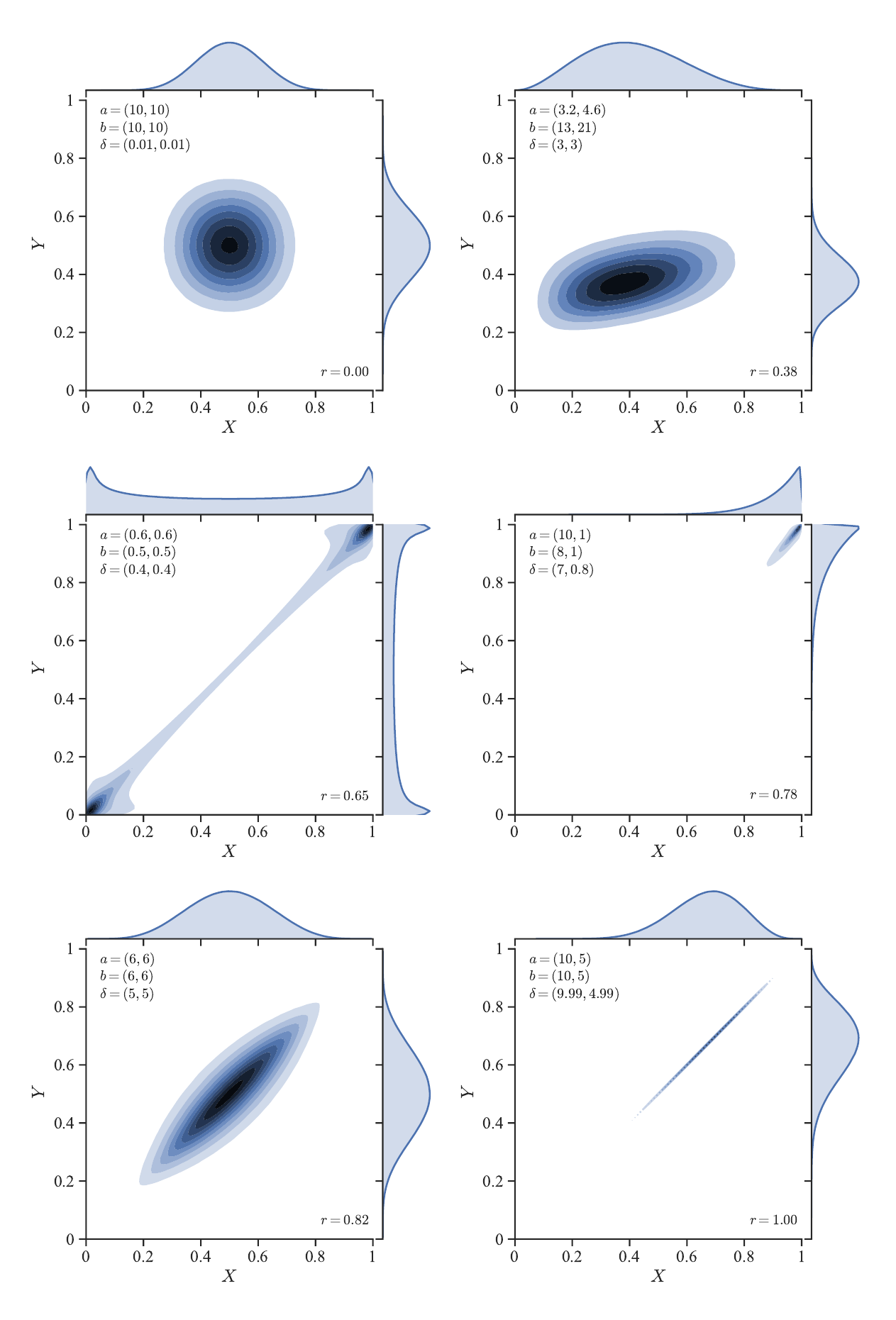}
	\caption{Joint densities of bivariate beta distributions with selected parameter values. The plots were created with kernel density estimation based on 10 million samples of the respective distributions. Note that the smoothing is inaccurate at the borders and produces artifacts close to zero and one as a consequence of smoothing with a symmetric kernel.}
	\label{fig:beta_densities}
\end{figure}

$X$ and $Y$ in (\ref{eq:uvw}) are linear transformations of $X^{\prime}$ and $Y^{\prime}$. Given $W_1,W_2,W_3$ it is easy to recover  $X^{\prime}$ and $Y^{\prime}$ from observed $X$ and $Y$,
\begin{align}
\begin{split}
X^{\prime} & =\frac{X-W_{3}}{W_{1}-W_{3}}=\frac{\left|X-W_{3}\right|}{\left|W_{1}-W_{3}\right|}=f_1(X)\\
Y^{\prime} & =\frac{Y-W_{3}}{W_{2}-W_{3}}=\frac{\left|Y-W_{3}\right|}{\left|W_{2}-W_{3}\right|}=f_2(Y).
\end{split}
\end{align}
As $X^{\prime}$ and $Y^{\prime}$ jointly follow the Olkin-Liu distribution \cite{olkin2003},
\begin{align}\label{eq:olkin_liu_joint}
p^{\prime}(x^{\prime},y^{\prime}) = \frac{\left(x^{\prime}\right)^{\upsilon_{1}-1}\left(1-x^{\prime}\right)^{\upsilon_{2}+\upsilon_{3}-1}\left(y^{\prime}\right)^{\upsilon_{2}-1}\left(1-y^{\prime}\right)^{\upsilon_{1}+\upsilon_{3}-1}}{B(\upsilon_{1},\upsilon_{2},\upsilon_{3})\left(1-x^{\prime}y^{\prime}\right)^{\upsilon_{1}+\upsilon_{2}+\upsilon_{3}}},
\end{align}
where $B(\upsilon_{1},\upsilon_{2},\upsilon_{3})=\frac{\Gamma(\upsilon_1)\Gamma(\upsilon_2)\Gamma(\upsilon_3)}{\Gamma(\upsilon_1+\upsilon_2+\upsilon_3)}$,
the joint distribution of $X$ and $Y$ given $W_1,W_2,W_3$ is
\begin{align} \label{eq:joint}
&p(x,y\mid w_1,w_2,w_3) =\left|\frac{df_1(x)}{dx}\frac{df_2(y)}{dy}\right|p^{\prime}(f_1(x),f_2(y)) \notag\\
&= \frac{\left|w_{1}-w_{3}\right|\left|w_{2}-w_{3}\right|}{B(\upsilon_{1},\upsilon_{2},\upsilon_{3})}\\
&\hspace{0.3cm}\cdot\frac{\left|x-w_{3}\right|^{\upsilon_{1}-1}\left|x-w_{1}\right|^{\upsilon_{2}+\upsilon_{3}-1}\left|y-w_{3}\right|^{\upsilon_{2}-1}\left|y-w_{2}\right|^{\upsilon_{1}+\upsilon_{3}-1}}{\left(\left|w_{1}-w_{3}\right|\left|w_{2}-w_{3}\right|-\left|x-w_{3}\right|\left|y-w_{3}\right|\right)^{\upsilon_{1}+\upsilon_{2}+\upsilon_{3}}} \notag
\end{align}
with $x$ between $w_1$ and  $w_3$ and $y$ between $w_2$ and $w_3$. We have not been able to integrate out $w_1,w_2,w_3$ from their joint distribution with $x$ and $y$. However, we suspect that even if the joint density for $X$ and $Y$ could be expressed in terms of special functions, computing those might not be efficient enough for parameter inference for which we will resort to moment matching. Example plots with smoothed samples for the joint density are shown in Figure \ref{fig:beta_densities} for several parameter settings showing different marginal distributions for $X$ and $Y$ and different correlations between $X$ and $Y$.

%%%%%%%%%%%%%%%%%%%%%%%%%%%%%%%%%%%%%%%%
\subsection{Moments} \label{sec:6_params_moments}
%%%%%%%%%%%%%%%%%%%%%%%%%%%%%%%%%%%%%%%%

As the marginal distributions for $X$ and $Y$ are beta-distributed, their moments are readily available, even in closed form. Computation of the product moments $\text{E}(X^kY^l)$ is more challenging but can be realized with help of the work of Olkin and Liu \cite{olkin2003}. Looking at construction (\ref{eq:uvw}), $X$ and $Y$ are a linear combination of independent beta-distributed random variables $W_1,W_2,W_3$ with weights $X^{\prime}$ and $Y^{\prime}$. Thus, we can express the product moments as
\begin{align} \label{eq:product_moments_xy_kl}
\text{E}(X^kY^l) = \text{E}\left(\left(X^{\prime}W_{1}+(1-X^{\prime})W_{3}\right)^k \left(Y^{\prime}W_{2}+(1-Y^{\prime})W_{3}\right)^l \right).
\end{align}
Since $X^{\prime}$ and $Y^{\prime}$ are independent of $W_1,W_2,W_3$, it is possible to compute the expectation if the moments of $W_1, W_2, W_3, X', Y'$ and the product moments of $X^{\prime}$ and $Y^{\prime}$ are known. $W_1, W_2, W_3$ as well as the marginals of $X'$ and $Y'$ are beta-distributed, so their moments can be computed straightforwardly in closed form. Furthermore, $X^{\prime}$ and $Y^{\prime}$ are jointly Olkin-Liu distributed according to (\ref{eq:olkin_liu_joint}), and Olkin and Liu \cite{olkin2003} have shown how to compute their product moments. However, note that the derivation of $\text{E}((X^{\prime})^k (Y^{\prime})^l)$ in equation (2.2) in \cite{olkin2003} is incorrect and should read
\begin{align}
\text{E}((X^{\prime})^k (Y^{\prime})^l) & =\sum_{j=0}^{\infty}dA(j)\frac{B(\upsilon_1+k+j,\upsilon_2+\upsilon_3)}{B(\upsilon_1+j,\upsilon_2+\upsilon_3)}\frac{B(\upsilon_2+l+j,\upsilon_1+\upsilon_3)}{B(\upsilon_2+j,\upsilon_1+\upsilon_3)} \label{eq:22_olkin}\\
%& =\sum_{j=0}^{\infty}dA(j)\frac{\Gamma(\upsilon_1+k+j)\Gamma(\upsilon_2+\upsilon_3)}{\Gamma(\Upsilon+k+j)} \notag \\
%&
%\hspace{1.1cm}\frac{\Gamma(\Upsilon+j)}{\Gamma(\upsilon_1+j)\Gamma(\upsilon_2+\upsilon_3)}\frac{\Gamma(\upsilon_2+l+j)\Gamma(\upsilon_1+\upsilon_3)}{\Gamma(\Upsilon+l+j)}\frac{\Gamma(\Upsilon+j)}{\Gamma(\upsilon_2+j)\Gamma(\upsilon_1+\upsilon_3)}\\
& =\frac{\Gamma(\upsilon_1+\upsilon_3)\Gamma(\upsilon_2+\upsilon_3)\Gamma(\upsilon_1+k)\Gamma(\upsilon_2+l)\Gamma(\Upsilon)}{\Gamma(\upsilon_1)\Gamma(\upsilon_2)\Gamma(\upsilon_3)\Gamma(\Upsilon+k)\Gamma(\Upsilon+l)} \notag \\
&\hspace{0.5cm}\sum_{j=0}^{\infty}\frac{(\upsilon_1+k)_{j}(\upsilon_2+l)_{j}(\Upsilon)_{j}}{(\Upsilon+k)_{j}(\Upsilon+l)_{j}}\frac{1}{j!} \label{eq:olkin_error1} \\
&= h \cdot  {}_3F_2(\upsilon_1+k,\upsilon_2+l,\Upsilon;\Upsilon+k,\Upsilon+l;1) \label{eq:product_moment},
\end{align}
with
\begin{align}
d & =\frac{\Gamma(\upsilon_1+\upsilon_3)\Gamma(\upsilon_2+\upsilon_3)}{\Gamma(\upsilon_3)\Gamma(\Upsilon)} \label{eq:d_olkin},\\
A(j) & =\frac{\Gamma(\upsilon_1+j)}{\Gamma(\upsilon_1)}\frac{\Gamma(\upsilon_2+j)}{\Gamma(\upsilon_2)}\frac{\Gamma(\Upsilon)}{\Gamma(\Upsilon+j)}\frac{1}{j!} \label{eq:Aj_olkin}\\
h & =\frac{\Gamma(\upsilon_1+\upsilon_3)\Gamma(\upsilon_2+\upsilon_3)\Gamma(\upsilon_1+k)\Gamma(\upsilon_2+l)\Gamma(\Upsilon)}{\Gamma(\upsilon_1)\Gamma(\upsilon_2)\Gamma(\upsilon_3)\Gamma(\Upsilon+k)\Gamma(\Upsilon+l)}\label{eq:olkin_error2}
\end{align}
where $\Upsilon=\upsilon_1+\upsilon_2+\upsilon_3$, ${}_pF_q$ is the generalized hypergeometric function. Equations (\ref{eq:22_olkin}), (\ref{eq:d_olkin}), and (\ref{eq:Aj_olkin}) are taken directly from \cite{olkin2003} with $a=\upsilon_1, b=\upsilon_2, c=\upsilon_3$. Equations (\ref{eq:olkin_error1}), (\ref{eq:product_moment}) and (\ref{eq:olkin_error2}) are our corrections of their equations.

%%%%%%%%%%%%%%%%%%%%%%%%%%%%%%%%%%%%%%%%
\subsection{Correlation and Covariance} \label{sec:correlation}
%%%%%%%%%%%%%%%%%%%%%%%%%%%%%%%%%%%%%%%%

The correlation $r$ between $X$ and $Y$ is
\begin{align} \label{eq:correlation}
r & =\frac{\text{Cov}(X,Y)}{\sqrt{\text{Var}(X)\text{Var}(Y)}}
\end{align}
with the known variances of the beta marginals
\begin{align}
\begin{split}
\text{Var}(X) & =\frac{a_1a_2}{(a_1+a_2)^{2}(a_1+a_2+1)}
=\frac{(\alpha_{1}+\delta_{1})(\alpha_{2}+\delta_{2})}{(\alpha_{1}+\delta_{1}+\alpha_{2}+\delta_{2})^{2}(\alpha_{1}+\delta_{1}+\alpha_{2}+\delta_{2}+1)}\\
\text{Var}(Y) &
=\frac{b_1b_2}{(b_1+b_2)^{2}(b_1+b_2+1)} =\frac{(\beta_{1}+\delta_{1})(\beta_{2}+\delta_{2})}{(\beta_{1}+\delta_{1}+\beta_{2}+\delta_{2})^{2}(\beta_{1}+\delta_{1}+\beta_{2}+\delta_{2}+1)}.
\end{split}
\end{align} 

For the covariance of $X$ and $Y$ Magnussen \cite{magnussen2004} gives an approximate solution, namely
\begin{align} \label{eq:magnussen_cov}
\text{Cov}(X,Y) &\approx \frac{a_1 a_2 \delta_2 + (1+b_1)(1+b_2)\delta_1}{(a_1+b_1)(a_2+b_2)(1+a_1+b_1)(1+a_2+b_2)},
\end{align}
where $a_1$, $a_2$, $b_1$, and $b_2$ are defined based on  $\alpha_1, \alpha_2, \beta_1, \beta_2, \delta_1, \delta_2$ as in (\ref{eq:ab}). This approximation is inaccurate for small values of these parameters, e.g. for $a_1 = a_2 = b_1 = b_2 = 4, \delta_1 = \delta_2 = 3$, the approximated covariance is Cov$(X,Y) = 0.024$, while the true covariance computed from $10^6$ samples of the bivariate beta distribution is Cov$(X,Y) = 0.020$. This difference results in an overestimated correlation of $r = 0.854$ as opposed to the true correlation of $r = 0.730$. Even more worryingly, for $a_1 = a_2 = b_1 = b_2 = 1, \delta_1 = \delta_2 = \frac{4}{5}$, the approximated covariance is Cov$(X,Y) = \frac{1}{9}$, which results in an estimated correlation of $r=\frac{4}{3}$, which is greater than 1 and therefore wrong by definition.

Given the connection to the Olkin-Liu distribution \cite{olkin2003}, which we derived in Section \ref{sec:6_params_moments}, we therefore proceed to compute the exact covariance between $X$ and $Y$:
\begin{align}
\text{Cov}(X,Y) & = \text{E}(XY) - \text{E}(X)\text{E}(Y)
\end{align}
where 
\begin{align}
\begin{split}
\text{E}(X) &
=\frac{a_1}{a_1+a_2}
=\frac{\alpha_{1}+\delta_{1}}{\alpha_{1}+\delta_{1}+\alpha_{2}+\delta_{2}} \\
\text{E}(Y) &
=\frac{b_1}{b_1+b_2}
=\frac{\beta_{1}+\delta_{1}}{\beta_{1}+\delta_{1}+\beta_{2}+\delta_{2}}.
\end{split}
\end{align}
are readily available as the means of the beta marginals. We can compute $\text{E}(XY)$ from (\ref{eq:product_moments_xy_kl}) with $k=l=1$, which results in
\begin{align}
\begin{split}
\text{E}(XY) & =\text{E}((X'W_{1}+(1-X')W_{3})(Y'W_{2}+(1-Y')W_{3}))\\
& =\text{E}(X'Y')\text{E}(W_{1})\text{E}(W_{2}) + \left(\text{E}(X')-\text{E}(X'Y')\right)\text{E}(W_{1})\text{E}(W_{3}) + \\
&\hspace{0.4cm}\left((\text{E}(Y')-\text{E}(X'Y')\right)\text{E}(W_{2})\text{E}(W_{3}) + \left(1-\text{E}(Y')-\text{E}(X')+\text{E}(X'Y')\right)\text{E}(W_{3}^{2})\\
\end{split}
\end{align}
with the moments of the beta marginals from (\ref{eq:uw}) 
\begin{align}
\begin{split}
\text{E}(W_{1}) & =\frac{\alpha_{1}}{\alpha_{1}+\alpha_{2}}\\
\text{E}(W_{2}) & =\frac{\beta_{1}}{\beta_{1}+\beta_{2}}\\
\text{E}(W_{3}) & =\frac{\delta_{1}}{\delta_{1}+\delta_{2}}\\
\text{E}(W_{3}^{2}) & =\frac{\delta_{1}(\delta_{1}+1)}{(\delta_{1}+\delta_{2}+1)(\delta_{1}+\delta_{2})}\\
\end{split}
\end{align}
and (\ref{eq:uvw})
\begin{align}
\begin{split}
\text{E}(X') & =\frac{\upsilon_1}{\upsilon_1+\upsilon_3}\\
\text{E}(Y') & =\frac{\upsilon_2}{\upsilon_2+\upsilon_3}
\end{split}
\end{align}
with $\upsilon_1$, $\upsilon_2$, and $\upsilon_3$ as defined in (\ref{eq:upsilon}). $\text{E}(X'Y')$ can be specialized from the moments (\ref{eq:product_moment}) above with $k=l=1$:
\begin{align}
\begin{split}
\text{E}(X'Y') & =h\hspace{0.1cm}_{3}F_{2}(\upsilon_1+1,\upsilon_2+1,\Upsilon;\Upsilon+1,\Upsilon+1;1)\qquad\text{with} \label{eq:3f2_cor}\\
h & =\frac{\Gamma(\upsilon_1+\upsilon_3)\Gamma(\upsilon_2+\upsilon_3)}{\Gamma(\upsilon_3)\Gamma(\Upsilon+1)}\cdot
	\frac{\upsilon_1 \upsilon_2}{\Upsilon}
\end{split}
\end{align}
and $\Upsilon=\upsilon_1+\upsilon_2+\upsilon_3$ as before. According to \cite{olkin2003} there is no closed-form solution for (\ref{eq:product_moment}) and (\ref{eq:3f2_cor}). However, using the generalized hypergeometric function the product moments and thus the covariance between $X$ and $Y$ can easily be computed numerically (e.g. by using the hyper function in sympy or the HypergeometricPFQ function in Mathematica).

%%%%%%%%%%%%%%%%%%%%%%%%%%%%%%%%%%%%%%%%
\subsection{Parameter Inference} \label{sec:parameter_inference}
%%%%%%%%%%%%%%%%%%%%%%%%%%%%%%%%%%%%%%%%

\begin{figure}[b!]
	\includegraphics[width = 0.95\textwidth]{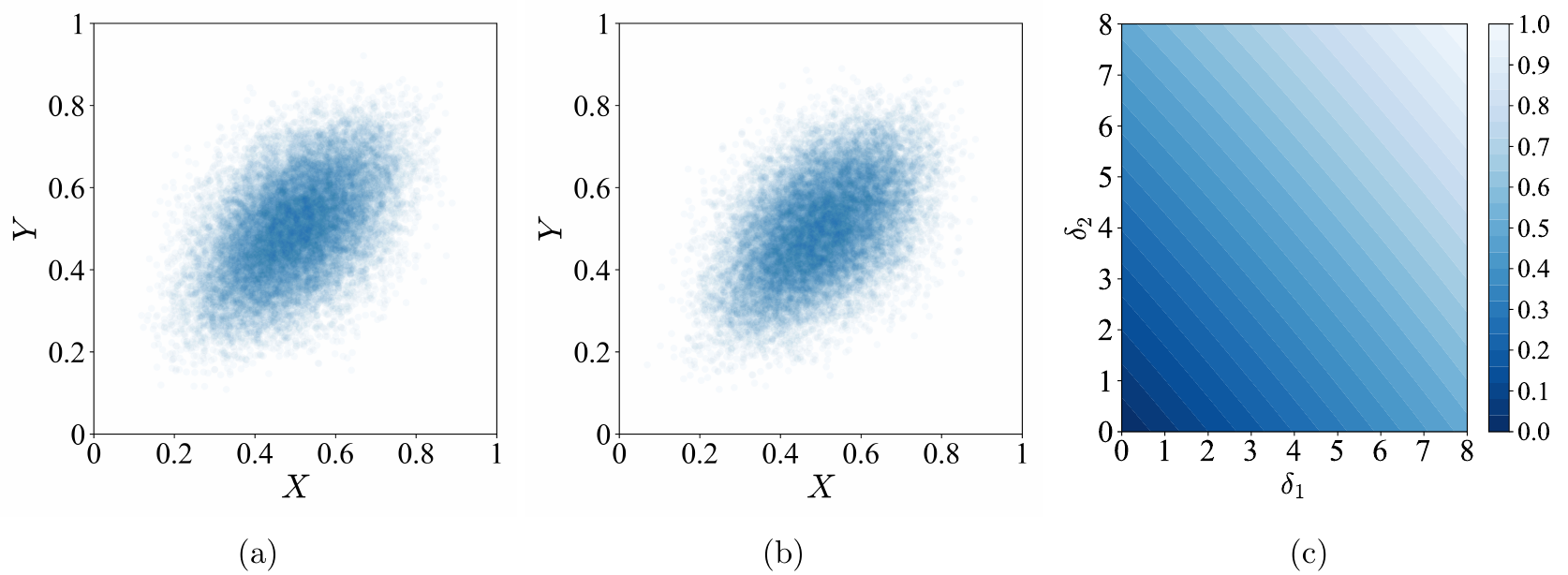}
	\caption{Different parameters can generate similar data. (a) and (b) show samples generated from two bivariate beta distributions with the same marginal parameters $a=b=(8,8)$, but different correlation parameters $\delta = (4,4)$ for (a) and $\delta = (6,2)$ for (b). Both, the data in (a) and (b) show a correlation of $r = 0.48$. Correspondingly, (c) shows the correlations generated by different combinations of $\delta_1$ and $\delta_2$ given marginal parameters $a = b = (8,8)$. Different combinations of $\delta_1$ and $\delta_2$ can lead to very similar correlations.}
	\label{fig:unidentifiable_example}
\end{figure}

Magnussen \cite{magnussen2004} used the method of moments to infer the parameters $a_1, a_2, b_1, b_2$ for the marginal distributions. Given these parameters he then matched the empirical correlation to the correlation for the parameters $\delta_1, \delta_2$ (given marginal parameters $a_1, a_2, b_1, b_2$) using the approximate solution for the covariance given in (\ref{eq:magnussen_cov}). However, as shown in Section \ref{sec:correlation} this approximation can be very inaccurate. An additional problem with the distribution is that very similar data can be generated with different parameter values, as one can see in Figure \ref{fig:unidentifiable_example}(a) and (b). Increasing $\delta_1$ and simultaneously decreasing $\delta_2$ or vice versa while keeping the marginal parameters $a_1, a_2$ and $b_1, b_2$ fixed, can result in very similar correlations, which is shown in Figure \ref{fig:unidentifiable_example}(c). Since two distributions with very different parameters can lead to extremely similar data, it is hard to statistically infer the parameters, $\delta_1$ and $\delta_2$, from data: The empirical correlation alone does not provide enough constraints and differences in higher moments can be subtle.

Therefore, in order to enable unambiguous parameter inference, we decided to constrain the 6-parameter bivariate beta distribution to five parameters: two for each marginal and one parameter to control the correlation. A reasonable way to constrain $\delta_1$ and $\delta_2$ is to set
\begin{align} \label{eq:constraint}
	\delta_2 &= \frac{\delta^{\text{max}}_2}{\delta^{\text{max}}_1} \delta_1
\end{align}
with $\delta^{\text{max}}_1 = \text{min}(a_1, b_1), \delta^{\text{max}}_2 = \text{min}(a_2, b_2)$, because this enables the maximum possible correlation between $X$ and $Y$ when the maximum values for $\delta_1$ and $\delta_2$ are attained and the shared component between $X$ and $Y$ is as big as it can be without violating the marginal constraints.

Using this constraint (\ref{eq:constraint}), the model-inherent constraints $\delta_1 < \delta^{\text{max}}_1, \delta_2 < \delta^{\text{max}}_2$, and the formula for the correlation derived in Section \ref{sec:correlation}, we can now optimize the parameters numerically to match the empirical moments. First, the marginal parameters $a_1, a_2, b_1, b_2$ are obtained using the standard procedure of moment matching for the beta distribution. Given the estimated marginal parameters, an estimate of $\delta_1$ (and with it $\delta_2$) can be obtained by minimizing the quadratic deviation between the numerically computed correlation and the empirical correlation.
\begin{figure}[b!]
	\centering
	\includegraphics[width = 0.55\textwidth]{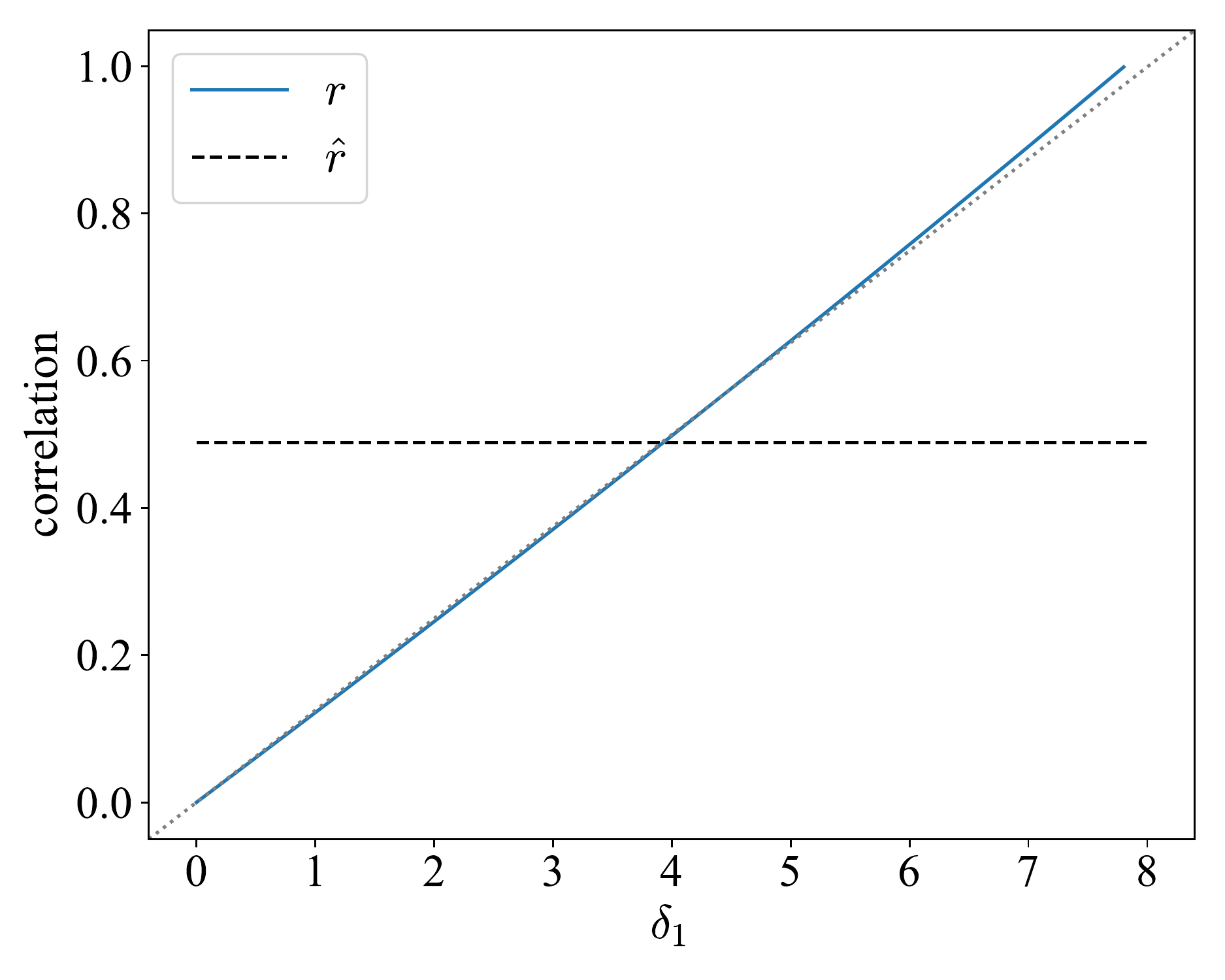}
	\caption{The correlations $r$ implied by different values of correlation parameter $\delta_1$ compared to the desired correlation $\hat{r} = 0.48$ given estimates of the marginal parameters $a_1 = 7.818, a_2 = 7.849, b_1 = 7.802, b_2 = 7.865$. The difference between $r$ and $\hat{r}$ is minimal for $\delta_1 = 3.93$. $\delta_2$ is not displayed since it can be computed from $\delta_1$ using constraint (\ref{eq:constraint}).}
	\label{fig:mm_error}
\end{figure}

As an example we used this numerical moment matching approach on 5000 data points generated with parameters $a_1 = a_2 = b_1 = b_2 = 8, \delta_1 = \delta_2 = 4$, equivalent to the data shown in Figure \ref{fig:unidentifiable_example}(a). We inferred the parameter values $a_1 = 7.818, a_2 = 7.849, b_1 = 7.802, b_2 = 7.865, \delta_1 = 3.93, \delta_2 = 3.954$. Figure \ref{fig:mm_error} shows the correlations implied by different values for parameter $\delta_1$ compared to the desired correlation of $\hat{r}=0.48$. As one can see, the difference between $r$ and $\hat{r}$ is minimal for the inferred $\delta_1 = 3.93$. Due to our constraint (\ref{eq:constraint}), $\delta_2 = 3.954$ can be computed from $\delta_1$. Also, note that in this case there is an almost linear relationship between $\delta_1$ and $r$.

%%%%%%%%%%%%%%%%%%%%%%%%%%%%%%%%%%%%%%%%%%%%%%%%%%%%%%%%%%%%%%%%%%%%%%%%%%%%%%%%
\section{Generalizations}
%%%%%%%%%%%%%%%%%%%%%%%%%%%%%%%%%%%%%%%%%%%%%%%%%%%%%%%%%%%%%%%%%%%%%%%%%%%%%%%%

%%%%%%%%%%%%%%%%%%%%%%%%%%%%%%%%%%%%%%%%
\subsection{Multivariate Beta Distribution}
%%%%%%%%%%%%%%%%%%%%%%%%%%%%%%%%%%%%%%%%

The bivariate beta distribution can be generalized to a multivariate beta distribution with more than two correlated beta-distributed random variables. We show this exemplarily for three dimensions. From 12 independent gamma-distributed random variables $A_1, A_2, B_1, B_2, C_1, C_2, D_1, D_2, E_1, E_2, F_1, F_2$ with distributions\\
\begin{minipage}{0.5\textwidth}
	\begin{align}
	A_i &\sim \text{Gamma}(\alpha_i,1)\quad i = 1,2 \notag \\
	B_i &\sim \text{Gamma}(\beta_i,1)\hspace{0.02cm}\quad i = 1,2 \notag \\
	C_i &\sim \text{Gamma}(\gamma_i,1)\hspace{0.03cm}\quad i = 1,2 \notag
	\end{align}
\end{minipage}
\begin{minipage}{0.5\textwidth}
	\begin{align}
	D_i &\sim \text{Gamma}(\delta_i,1)\quad i = 1,2 \notag \\
	E_i &\sim \text{Gamma}(\epsilon_i,1)\quad i = 1,2 \\
	F_i &\sim \text{Gamma}(\phi_i,1)\hspace{0.3cm} i = 1,2 \notag
	\end{align}
\end{minipage}\\\\
we construct the three multivariate-beta-distributed random variables
\begin{align}\label{eq:multivariate_beta}
X &= \frac{A_1 + D_1 + F_1}{A_1 + A_2 + D_1 + D_2 + F_1 + F_2}, \notag \\
Y &= \frac{B_1 + D_1 + E_1}{B_1 + B_2 + D_1 + D_2 + E_1 + E_2}, \\
Z &= \frac{C_1 + E_1 + F_1}{C_1 + C_2 + E_1 + E_2 + F_1 + F_2}. \notag
\end{align}
The marginal distributions of $X$, $Y$ and $Z$ are Beta($X$; $\alpha_1 + \delta_1 + \phi_1$, $\alpha_2 + \delta_2 + \phi_2$), Beta($Y$; $\beta_1 + \delta_1 + \epsilon_1$, $\beta_2 + \delta_2 + \epsilon_2$) and Beta($Z$; $\gamma_1 + \epsilon_1 + \phi_1$, $\gamma_2 + \epsilon_2 + \phi_2$). Correlations between $X$ and $Y$ are obtained through the common parameters $\delta_1$ and $\delta_2$, between $Y$ and $Z$ through parameters $\epsilon_1$ and $\epsilon_2$ and between $X$ and $Z$ through parameters $\phi_1$ and $\phi_2$. Generalizing the bivariate beta distribution to a $k$-dimensional multivariate beta distribution requires $2(k + \binom{k}{2})$ parameters.

%%%%%%%%%%%%%%%%%%%%%%%%%%%%%%%%%%%%%%%%
\subsection{Correlated Dirichlet Distribution}
%%%%%%%%%%%%%%%%%%%%%%%%%%%%%%%%%%%%%%%%

\begin{figure}
	\includegraphics[width=0.95\textwidth]{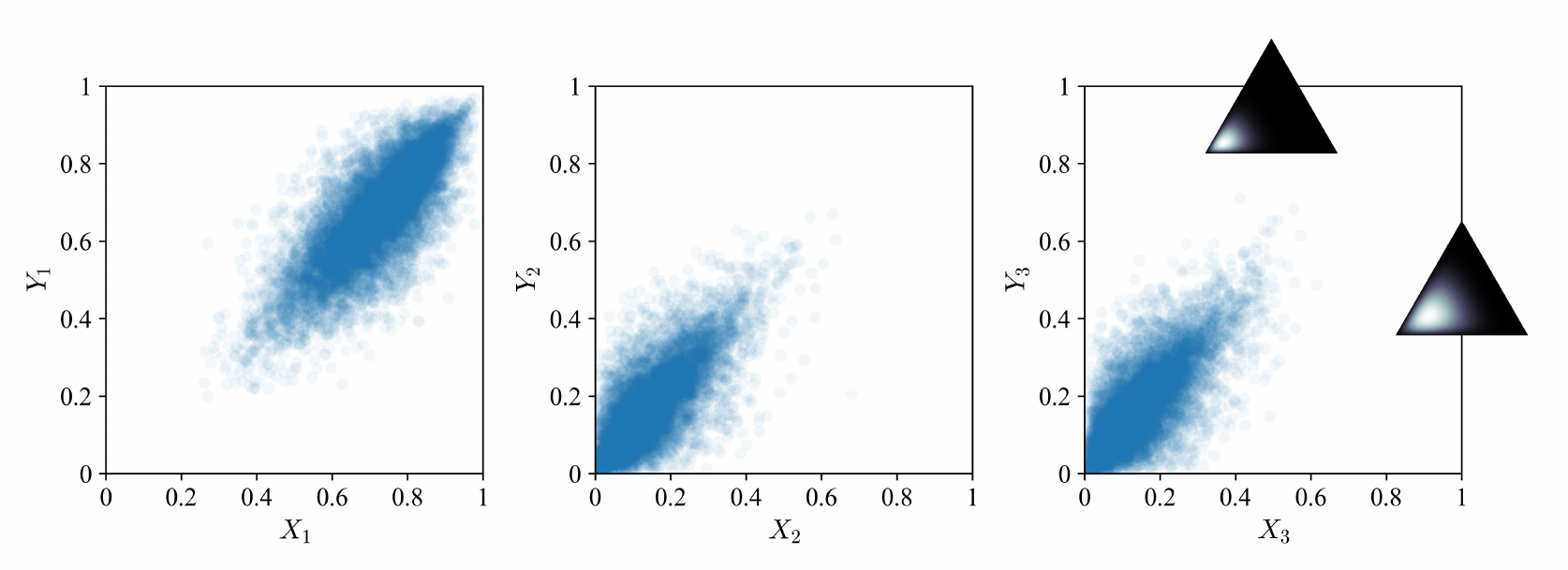}
	\caption{Data samples generated from a correlated Dirichlet distribution with parameters $\alpha = (2.5,0.5,0.5), \beta = (0.5,0.5,0.5), \delta = (7.5, 1.5, 1.5)$. The marginal distributions, displayed by the simplexes, are Dirichlet($X; 10,2,2$) and Dirichlet($Y; 8,2,2$). Correlations between $X$ and $Y$ in the respective dimensions are $r_1 = 0.77, r_2 = 0.76, r_3 = 0.76$.}
	\label{fig:cor_dirichlet_example}
\end{figure}

The bivariate beta distribution described in Section \ref{sec:6_params} models two correlated random variables $X$ and $Y$ with beta-distributed marginal distributions. We now generalize this case to two correlated random vectors $\bm{X}=(X_1, \dots, X_k)$ and $\bm{Y}=(Y_1, \dots, Y_k)$ with the two marginals being Dirichlet-distributed, resulting in a correlated Dirichlet distribution. First note that the beta distribution can be seen as a two-dimensional distribution over two values in $[0,1]$ which sum up to 1. For the bivariate beta distribution in (\ref{eq:bivariate_beta_6}) $X_1$ and $Y_1$ as well as $X_2=1-X_1$ and $Y_2=1-Y_1$ are generated from six gamma-distributed random variables,
\begin{align}
\begin{split}
X_1 &= \frac{A_1 + D_1}{A_1 + A_2 + D_1 + D_2} \quad\mathrm{and}\quad X_2 = 1-X_1 = \frac{A_2 + D_2}{A_1 + A_2 + D_1 + D_2} \\
Y_1 &= \frac{B_1 + D_1}{B_1 + B_2 + D_1 + D_2} \quad\mathrm{and}\quad Y_2 = 1-Y_1 = \frac{B_2 + D_2}{B_1 + B_2 + D_1 + D_2}.
\end{split}	
\end{align}
From this viewpoint, it is straightforward to generalize the bivariate beta distribution to a correlated Dirichlet distribution. A $k$-dimensional correlated Dirichlet distribution can be constructed from $3k$ gamma-distributed random variables $A_1, \dots, A_k, B_1, \dots, B_k, D_1, \dots, D_k$ with $3k$ parameters $\alpha_1, \dots, \alpha_k, \beta_1, \dots \beta_k, \delta_1, \dots \delta_k$ distributed according to
\begin{align}
\begin{split}
A_i &\sim \text{Gamma}(\alpha_i, 1)\quad i = 1,\dots,k \\
B_i &\sim \text{Gamma}(\beta_i, 1)\hspace{0.02cm}\quad i = 1,\dots,k \\
D_i &\sim \text{Gamma}(\delta_i, 1)\quad\hspace{0.06cm} i = 1,\dots,k.
\end{split}
\end{align}
These random variables are used to construct the correlated Dirichlet-distributed random variables $\bm{X}=(X_1,\dots,X_k)$ and $\bm{Y}=(Y_1,\dots,Y_k)$ with\\
\begin{minipage}{0.49 \linewidth}
	\begin{align}
	X_i &= \frac{A_i + D_i}{\sum_{i=1}^{k}{A_i} + \sum_{i=1}^{k}{D_i}} \quad\qquad\mathrm{and}\notag	
	\end{align}
\end{minipage}
\begin{minipage}{0.49 \linewidth}
	\begin{align} \label{eq:correlated_dirichlet}
	Y_i &= \frac{B_i + D_i}{\sum_{i=1}^{k}{B_i} + \sum_{i=1}^{k}{D_i}}.
	\end{align}
\end{minipage}\\\\
The two resulting marginal distributions are Dirichlet($\bm{X}; \alpha_1 + \delta_1,  \dots, \alpha_k + \delta_k$) and Dirichlet($\bm{Y}; \beta_1 + \delta_1, \dots, \beta_k + \delta_k$). An example of data generated from a 3-dimensional correlated Dirichlet distribution can be seen in Figure \ref{fig:cor_dirichlet_example}.
%

%%%%%%%%%%%%%%%%%%%%%%%%%%%%%%%%%%%%%%%%%%%%%%%%%%%%%%%%%%%%%%%%%%%%%%%%%%%%%%%%
\section*{Funding}
%%%%%%%%%%%%%%%%%%%%%%%%%%%%%%%%%%%%%%%%%%%%%%%%%%%%%%%%%%%%%%%%%%%%%%%%%%%%%%%%

This work has been supported by the German Federal Ministry of Education and Research (BMBF) in the project KoBo34 (project no. 16SV7984).

%%%%%%%%%%%%%%%%%%%%%%%%%%%%%%%%%%%%%%%%%%%%%%%%%%%%%%%%%%%%%%%%%%%%%%%%%%%%%%%%

\bibliographystyle{tfnlm}
\bibliography{references_trick_jaekel_rothkopf}

\end{document}